\documentclass[preprint,12pt]{elsarticle}

\usepackage{amssymb}

\journal{Journal of Number Theory}
\usepackage{hyperref} 
\usepackage{amsmath} 
\usepackage{mathrsfs}  
\usepackage{amsthm}
\usepackage{color}
\newtheorem{thm}{Theorem} 
\newtheorem{lem}[thm]{Lemma} 
\newtheorem{conj}[thm]{Conjecture} 
\newdefinition{rmk}{Remark} 
\newproof{pf}{Proof} 
\numberwithin{equation}{section}
\begin{document}

\begin{frontmatter}



\title{Siegel Zeros and the Hardy-Littlewood Conjecture}


\author{Yunan Wang}

\affiliation{organization={Clemson University},
            addressline={}, 
            city={Clemson},
            postcode={29634}, 
            state={South Carolina},
            country={United States}}
\begin{abstract}

In 2016, Fei \cite{fei2016application} established a bound on the Siegel zeros for real primitive Dirichlet characters modulo $q$, assuming the weak Hardy-Littlewood conjecture. Building on Fei's work, Jia \cite{jia2022conditional} demonstrated the same bound using a stronger version of the Hardy-Littlewood conjecture. In this paper, we present a slightly simplified approach to reprove their results.

\end{abstract}




\end{frontmatter}


\section{Introduction} \label{sec:intro}
Zero-free regions for Dirichlet L-functions is an important problem in number theory.  Let \(\chi\) be a complex character modulo \(q\), and let \(s = \sigma + \mathrm{i}t\). It is well-established that there exists a constant \(c_1 > 0\) such that the region defined by
\[
\sigma > 1 - \frac{c_1}{\log(q(|t| + 2))}
\]
is free of zeros for the Dirichlet \(L\)-function \(L(s, \chi)\). However, if \(\chi\) is a real character modulo \(q\), then there may be at most one simple real zero. If such a zero \(\beta\) exists, it is termed an exceptional zero or Siegel zero, and the character \(\chi\) is referred to as an exceptional character. In the case, both $\chi$ and $s$ are real much less is known. In 1935, Siegel \cite{page1935number} showed that for each $\varepsilon>0$, there is an ineffective constant $c_2(\varepsilon)>0$ such that for any real character $\chi \bmod q(q \geq 3)$, if $\beta$ is a real zero of $L(s, \chi)$, then
$$
\frac{c_2(\varepsilon)}{q^{\varepsilon}} \leq 1-\beta.
$$

Recently, studies have shown that assumptions related to the Goldbach conjecture can yield significant bounds for exceptional zeros.  In 1923, Hardy and Littlewood \cite{hardy1923some} conjectured that for the function \(R(n)\), which represents the number of ways \(n\) can be expressed as a sum of two primes, the following asymptotic formula holds:
\[
R(n) \sim 2C_{2} \left(\prod_{\substack{p \mid n \\ p \geq 3}} \frac{p-1}{p-2}\right) \frac{n}{\log^2(n)},
\]
where \(C_{2}\) is the twin prime constant, defined by 
\begin{equation}
    C_{2} := \prod_{\substack{\text{prime} \\ p \geq 3}}\left(1-\frac{1}{(p-1)^{2}}\right) \approx 0.660161815846869573 \dots. \nonumber
\end{equation}

A weakened version of the Hardy-Littlewood conjecture can be described as follows:
\begin{conj} \label{conj1}
    There exists an absolute constant $c_3>0$ such that for all the even integer $n \geq 4$, we have
$$
R(n) \geq \frac{c_3 n}{\log ^2 n}.
$$
\end{conj}
In 2016, Fei\cite{fei2016application} showed that if conjecture \ref{conj1} holds, then the bounds will be improved by the following theorem. 
\begin{thm}[Fei(2016)] \label{fei's them2}
   Let $q$ be a prime number congruent to 3 modulo 4 and let $\chi$ be the real primitive character of modulus $q$. Suppose that the Dirichlet L-function $L(s, \chi)$ has a real exceptional zero $\beta$. If conjecture \ref{conj1}  is correct then there is a number $c_4>0$ such that for any such $q$ and $\beta$, we have
$$
\beta \leq 1-\frac{c_4}{\log ^2 q}
$$
\end{thm}

In 2022, Using the similar idea of Fei, under an assumption of a little stronger conjecture of Conjecture \ref{conj1}, Jia \cite{jia2022conditional} showed the same bound. 
\begin{conj}\label{conj-jia1}
 Suppose that $x$ is sufficiently large, $q \leq \frac{x}{4}$. There is an absolute constant $c_5>0$ such that for the even integers $n\left(\frac{x}{2}<n \leq x, q \mid n\right)$, we have
$$
R(n) \geq \frac{c_5 n}{\varphi(n)} \cdot \frac{n}{\log ^2 n}
$$
\end{conj}
\begin{thm}[Jia(2022)] \label{Jia's them4}
    Suppose that Conjecture \ref{conj-jia1} holds true. Let $q$ be a sufficiently large composite number, $\chi$ be the real primitive character $\bmod q$ with $\chi(-1)=-1$ and $\beta$ be the Siegel zero of $L(s, \chi)$. Then there is an effective absolute constant $c_6>0$ such that
$$
\beta \leq 1- \frac{c_6}{\log ^2 q} . 
$$
\end{thm}
In the same year 2022, followed by Fei's work, Friedlander, Goldston, Iwaniec and  Suriajaya \cite{friedlander2022exceptional} generalized Fei and Jia's result, showed that the assumption of a weak form of the HardyLittlewood conjecture on the Goldbach problem suffices to disprove the possible existence of exceptional zeros of Dirichlet $L$-functions.

Define $$G(n)=\sum_{\substack{m_1+m_2=n \\ 2 \nmid m_1 m_2}} \Lambda\left(m_1\right) \Lambda\left(m_2\right) \text{and} $$

$$\mathfrak{S}(n)=2 C_2\prod_{\substack{p \mid n \\ p>2}}\left(1+\frac{1}{(p-2)}\right) .$$

\begin{conj}\label{conj-gpy}
    For all sufficiently large even $n$ we have
    $$
\delta \mathfrak{S}(n) n<G(n)<(2-\delta) \mathfrak{S}(n) n
$$
for some fixed $0<\delta<1$.
\end{conj}

\begin{thm}[Friedlander, Goldston, Iwaniec, and Surajaya(2022)]
    Assume that the Weak Hardy-Littlewood-Goldbach Conjecture \ref{conj-gpy} holds for all sufficiently large even $n$. Then there are no zeros of any Dirichlet L-function in the region $\sigma \geqslant 1-c / \log q(|t|+2)$ with a positive constant $c$ which is now allowed to depend on $\delta$.
\end{thm}
Fei derived Theorem \ref{fei's them2} by comparing two different results of $\sum_{k=1}^q\left(\sum_{3 \leq p \leq x} e\left(\frac{k p}{q}\right)\right)^2$, then use the weak Hardy-Littlewood conjecture.  Jia derived Theorem \ref{Jia's them4} using similar idea with slightly different weak Hardy-Littlewood conjecture. In this paper, we'll use a simpler method to get a similar result of Theorem \ref{fei's them2} and \ref{Jia's them4}.

\section{Some Lemmas}
\begin{lem}{(Prime number theorem for arithmetic progressions).} \label{lem7}
If there exists an exceptional character $\chi_1$ modulo $q$ and a real zero $\beta>1-\frac{c_{11}}{\log q}$ of $L(s, \chi)$ (where $c_{11}>0$ is a fixed positive constant independent of $q$ ), then for a constant $c_{10}>0$, such that $q \leq \exp \left(2 c_{10} \sqrt{\log x}\right)$, then if $(a,q)=1$, we have
$$\pi(x ; q, a)=\frac{\operatorname{li}(x)}{\varphi(q)}-\frac{\chi_1(a) \operatorname{li}\left(x^{\beta}\right)}{\varphi(q)}+O\left(x \exp \left(-c_{10} \sqrt{\log x}\right)\right)$$
\end{lem}

One could see [\cite{montgomery2007multiplicative}, p. 381, Corollary 11.20]


Dirichlet and Legendre conjectured, and de la Vallée Poussin proved the following lemma:
\begin{lem}\label{lem8}
If $(a,d)=1$, then
    $$\pi_{d, a}(x) \sim \frac{\operatorname{li}(x)}{\varphi(d)}$$
\end{lem}
This lemma shows if $(a_1,d)=1$ and $(a_2,d)=1$, then $\pi_{d, a_1}(x)$ and $\pi_{d, a_2}(x)$ are about the same size. 
\section{proof of theorem 2}

Let $q$ be a given positive prime number. We consider the sum 
\begin{equation}
    r_q(x)=\sum_{\substack{3 \leq p_1, p_2 \leq x \\ p_1+p_2 \equiv 0(\bmod q)}} 1 
\end{equation}
 in two different ways. 
\subsection{First estimation of \texorpdfstring{\( {r_q}(x) \)}{srq(x)}}

 We take $x=\exp \left(\frac{36}{c_{10}^2} \log ^2 q\right)$, where $c_{10}$ is defined in Lemma \ref{lem7}, since $q$ is sufficiently large, so is $x$. And we have $q=\exp \left(\frac{c_{10}}{6} \sqrt{\log x}\right)<\frac{x}{4}$.

 We first consider the lower bound of $r_q(x)$.
\begin{equation}
    r_q(x)=\sum_{\substack{3 \leq p_1, p_2 \leq x \\ p_1+p_2 \equiv 0(\bmod q)}} 1 \geq \sum_{n=1}^{\left[\frac{x}{2 q}\right]} \sum_{\substack{3 \leq p_1, p_2 \leq x \\ p_1+p_2=2 n q}} 1
\end{equation}
Under the assumption of Conjecture \ref{conj1}, we have:
\begin{equation}
    \sum_{n=1}^{\left[\frac{x}{2 q}\right]} \sum_{\substack{3 \leq p_1, p_2 \leq x \\ p_1+p_2=2 n q}} 1 \geq  \sum_{n=1}^{\left[\frac{x}{2 q}\right]}\frac{\delta 2 n q}{\log ^2 2 n q}\geq \frac{\delta x^2}{4 q\log ^2 x}+O\left(\frac{x}{\log ^2 x}\right)
\end{equation}
Since the error term comes from a positive term, so we have the lower bound:
\begin{equation} \label{di1ge}
    r_q(x)\geq \frac{\delta x^2}{4 q\log ^2 x}
\end{equation}
\subsection{Second estimation of \texorpdfstring{\( r_q(x) \)}{rq(x)}}

$\begin{aligned} r_q(x) & =\sum_{\substack{3 \leq p_1, p_2 \leq x \\ p_1+p_2 \equiv 0(\bmod q)}} 1\leq \sum_{i=1}^{q}\left(\sum_{\substack{ p \leq x \\ p \equiv i(\bmod q)}}1\sum_{\substack{ p \leq x \\ p \equiv q-i(\bmod q)}}1\right)\\&\leq q \left( \sum_{\substack{ p \leq x \\ p \equiv 1(\bmod q)}}1\right)^2=q\left(\frac{\operatorname{li}(x)}{\varphi(q)}-\frac{\chi(1)}{\varphi(q)} \operatorname{li}\left(x^\beta\right)+O\left(x \exp \left(-c_{10} \sqrt{\log x}\right)\right) \right)^2\\&=q\left(\frac{\operatorname{li}^2(x)}{\varphi^2(q)}-2\frac{\operatorname{li}(x)\chi(1)\operatorname{li}\left(x^\beta\right)}{\varphi^2(q)}+O(f(x))\right) \\\end{aligned}$

where $f(x)$= \text{max} $\{\frac{x^{2\beta - 2}}{\log x},\frac{1}{\sqrt{x^{c_{10}}}}\}$. And the second line can be derived from

Note that
$$
\operatorname{li}(x)=\frac{x}{\log x}+O\left(\frac{x}{\log ^2 x}\right) 
\text{and} \varphi{(q)}=q-1
$$

Then we have\\ \\
$\begin{aligned}
    & q\left(\frac{\operatorname{li}^2(x)}{\varphi^2(q)}-2\frac{\operatorname{li}(x)\chi(1)\operatorname{li}\left(x^\beta\right)}{\varphi^2(q)}+O(f(x))\right) \\
    & =\frac{1}{q}\frac{x^2}{\log ^2 x}+O\left(\frac{x^2}{\log ^3 x}\right)
-\frac{2}{q}\left(\frac{x}{\log x}+O\left(\frac{x}{\log ^2 x}\right)\right)\left(\frac{x^\beta}{\beta \log x}+O\left(\frac{x^\beta}{\log ^2 x}\right)\right)+O(f(x))\end{aligned}$
\\ \\
After simply, we have:\\
\begin{equation} \label{di2ge}
    r_p(s)\leq \frac{x^2}{q\log ^2 x}-2\frac{x^{1+\beta}}{q\log ^2 x}+O\left(
\frac{x^2}{\log^3 x}\right).
\end{equation}
Then we combine \ref{di1ge} and \ref{di2ge} and have: 
$$\frac{\delta x^2}{4 q\log ^2 x}\leq \frac{x^2}{q\log ^2 x}-2\frac{x^{1+\beta}}{q\log ^2 x}+O\left(
\frac{x^2}{\log^3 x}\right)$$
which yields 
$$\beta-1\leq \frac{\log(\frac{1}{2}-\frac{c_{11}}{8})}{\log x}\leq \frac{c_{10}^2\log(\frac{1}{2}-\frac{c_{11}}{8})}{36\log^2 q} $$
which is 
$$\beta \leq 1- \frac{c_{10}}{\log^2 q}$$
where $c_{10}>0$ is an effective absolute constant. \\
Based on our assumption that $q$ is a large prime. So for the above case $\varphi{(q)}=q-1$. 

In 2022, Jia \cite{jia2022conditional} employed a more stronger version of the Hardy-Littlewood conjecture to establish a result similar to Fei's \cite{fei2016application}, but for the case where $q$ is composite.

In the next section, we will employ a proof similar to the one used in our previous proof, but our method will be simpler than Jia's approach as we demonstrate Theorem 4.

\section{Proof of Theorem 4}

$\begin{aligned}
    r_q(x) &\geq \sum_{n=1}^{\left[\frac{x}{2 q}\right]} \sum_{\substack{3 \leq p_1, p_2 \leq x \\ p_1+p_2=2 n q}} 1 \geq \sum_{n=1}^{\left[\frac{x}{2 q}\right]}\frac{c_8(2nq)^2}{\varphi(2nq)\log^2(2nq)} \geq \frac{4c_8q^2}{\varphi(q)\log^2x}\sum_{n=1}^{\left[\frac{x}{2 q}\right]}\frac{n^2}{\varphi(2n)} \\
    &\geq \frac{4c_8q^2}{\varphi(q)\log^2x}\sum_{n=1}^{\left[\frac{x}{2 q}\right]}n \geq \frac{c_8x^2}{\varphi(q)\log^2x} \\
\end{aligned}$
\\So we have 
\begin{equation} \label{di4ge}
    r_q(x)\geq \frac{c_8x^2}{\varphi(q)\log^2x}
\end{equation}
On the other hand, \\ \\
$\begin{aligned}
    r_q(x) & \leq q\left(\frac{\operatorname{li}^2(x)}{\varphi^2(q)}-2\frac{\operatorname{li}(x)\chi(1)\operatorname{li}\left(x^\beta\right)}{\varphi^2(q)}+O(f(x))\right) \\
    & \leq \frac{x^2q}{\varphi ^2(q)\log ^2 x}+O\left(\frac{x^2}{\log ^3 x}\right)
-\frac{2q}{\varphi^2(q)}\left(\frac{x}{\log x}+O\left(\frac{x}{\log ^2 x}\right)\right)\left(\frac{x^\beta}{\beta \log x}+O\left(\frac{x^\beta}{\log ^2 x}\right)\right)\\
&+O(f(x))\\
\end{aligned}$
After simply, we have:\\
\begin{equation} \label{di3ge}
    r_q(s)\leq \frac{qx^2}{\varphi^2(q)\log ^2 x}-2\frac{qx^{1+\beta}}{\varphi^2(q)\log ^2 x}+O\left(
\frac{x^2}{\log^3 x}\right).
\end{equation}

Then we combine \ref{di3ge} and \ref{di4ge}, 
\begin{equation}\label{four}
    \frac{c_8x^2}{\varphi(q)\log^2x}\leq \frac{qx^2}{\varphi^2(q)\log ^2 x}-2\frac{qx^{1+\beta}}{\varphi^2(q)\log ^2 x}+O\left(
\frac{x^2}{\log^3 x}\right)
\end{equation}
Note that 
$$\varphi(q)=q \prod_{p \mid q}\left(1-\frac{1}{p}\right)$$
and equation \ref{four} gives
$$\frac{qx^{\beta-1}}{\varphi(q)}\leq \frac{q}{2\varphi(q)}-\frac{c_8}{2}$$
Let $\Tilde{c_9}=\prod_{p \mid q\left(1-\frac{1}{p}\right)}$ for a fixed large composite $q$.

So we have
$$\beta-1 \leq \frac{\log \left(\frac{1}{2}-\frac{c_{8}\Tilde{c_9}}{8}\right)}{\log x} \leq \frac{c_{10}^2 \log \left(\frac{1}{2}-\frac{c_{8}\Tilde{c_9}}{8}\right)}{36 \log ^2 q}$$

So we have $$\beta \leq 1- \frac{c_{12}}{\log^2 q}.$$ 

\bibliographystyle{elsarticle-num} 
\bibliography{references} 

\begin{thebibliography}{1}
\expandafter\ifx\csname url\endcsname\relax
  \def\url#1{\texttt{#1}}\fi
\expandafter\ifx\csname urlprefix\endcsname\relax\def\urlprefix{URL }\fi
\expandafter\ifx\csname href\endcsname\relax
  \def\href#1#2{#2} \def\path#1{#1}\fi

\bibitem{fei2016application}
J.~Fei, An application of the hardy--littlewood conjecture, Journal of Number Theory 168 (2016) 39--44.

\bibitem{jia2022conditional}
C.~H. Jia, On the conditional bounds for siegel zeros, Acta Mathematica Sinica, English Series 38~(5) (2022) 869--876.

\bibitem{page1935number}
A.~Page, On the number of primes in an arithmetic progression, Proceedings of the London Mathematical Society 2~(1) (1935) 116--141.

\bibitem{hardy1923some}
G.~H. Hardy, J.~E. Littlewood, Some problems of ‘partitio numerorum’; iii: On the expression of a number as a sum of primes, Acta Mathematica 44~(1) (1923) 1--70.

\bibitem{friedlander2022exceptional}
J.~Friedlander, D.~Goldston, H.~Iwaniec, A.~Suriajaya, Exceptional zeros and the goldbach problem, Journal of Number Theory 233 (2022) 78--86.

\bibitem{montgomery2007multiplicative}
H.~L. Montgomery, et~al., Multiplicative Number Theory, no.~97 in Graduate Texts in Mathematics, Springer, 2007.

\end{thebibliography}

\end{document}